\documentclass[12pt, hidelinks]{article}
\usepackage[utf8]{inputenc}
\usepackage{amsmath, amssymb, url, float, graphicx, float}
\usepackage{algorithm, algpseudocode}
\usepackage{fullpage}
\usepackage[small,bf]{caption}
\usepackage{tikz}
\usetikzlibrary{patterns, positioning}

\newcommand{\eps}{\varepsilon}

\newcommand{\tmax}{\theta^\mathrm{max}}
\newcommand{\tmin}{\theta^\mathrm{min}}
\newcommand{\gmax}{g^\mathrm{max}}
\newcommand{\gmin}{g^\mathrm{min}}

\newcommand{\C}{\mathcal{C}}

\newcommand{\ones}{\mathbf 1}
\newcommand{\reals}{{\mbox{\bf R}}}

\newcommand{\diag}{\mathop{\bf diag}}

\newcommand{\argmin}{\mathop{\rm argmin}}

\newcommand{\sign}{\mathop{\bf sign}}

\newcommand{\eg}{{\it e.g.}}
\newcommand{\ie}{{\it i.e.}}

\newcommand{\BEAS}{\begin{eqnarray*}}
\newcommand{\EEAS}{\end{eqnarray*}}
\newcommand{\BEA}{\begin{eqnarray}}
\newcommand{\EEA}{\end{eqnarray}}
\newcommand{\BEQ}{\begin{equation}}
\newcommand{\EEQ}{\end{equation}}
\newcommand{\BIT}{\begin{itemize}}
\newcommand{\EIT}{\end{itemize}}

\title{A New Heuristic for Physical Design}
\date{February 2020}
\author{Guillermo Angeris \\
\texttt{\small angeris@stanford.edu}
\and
Jelena Vu\v{c}kovi\'c
\\
\texttt{\small jela@stanford.edu}
\and
Stephen Boyd \\
\texttt{\small boyd@stanford.edu}
}

\begin{document} 
\maketitle 
\begin{abstract} 
In a physical design problem, the designer chooses values of
some physical parameters, within limits, to optimize the resulting field.
We focus on the specific case in which each physical design parameter
is the ratio of two field variables.
This form occurs for photonic design with real scalar fields, diffusion-type
systems, and others.
We show that such problems can be reduced to a convex optimization
problem, and therefore efficiently solved globally, given the sign of 
an optimal field at every point.
This observation
suggests a heuristic, in which the signs of the field are iteratively updated.
This heuristic appears to have good practical performance on
diffusion-type problems (including thermal design and resistive circuit design)
and some control problems, while exhibiting moderate performance on photonic
design problems.  We also show in many practical
cases there exist globally optimal designs whose design parameters are
maximized or minimized at each point in the domain, \ie, that there is a
discrete globally optimal structure.  
\end{abstract}

\section{Introduction}
Computer-aided physical design has become an important tool in many fields including photonics~\cite{moleskyInverseDesignNanophotonics2018, suNanophotonicInverseDesign2019}, mechanical design~\cite{haftkaElementsStructuralOptimization2012}, circuit design~\cite{gielenComputeraidedDesignAnalog2000, liuAnalogCircuitOptimization2009}, and thermal design~\cite{dboukReviewEngineeringDesign2017, haertelDesignThermalSystems2018}. In many cases, the design problem is formulated as a constrained nonconvex optimization problem which is then approximately minimized using local optimization methods such as ADMM~\cite{luNanophotonicComputationalDesign2013}, evolutionary algorithms~\cite{liuAnalogCircuitOptimization2009}, and the method of moving asymptotes~\cite{dboukReviewEngineeringDesign2017}, among many others.

More generally, a physical design problem can be phrased in the following way:
we are allowed to choose some design parameters (\eg, the permittivity in
photonic design or the conductances in diffusion design) at each point in a
domain, within some limits, in order to minimize an objective function of the
field (this can be, \eg, the electric field in photonic design, or a vector
containing the potentials, flows, and potential differences in diffusion
design). The constraints specify the physics of the problem, connecting the
design variables to the field variables (\eg, Maxwell's equations in photonics,
or a diffusion equation such as the heat equation in diffusion design). We
note that, in many cases, the physics constraints are linear equations in the
field variables (when the design parameters are held constant), and linear equations in
the design parameters (when the fields are held constant), which has led to
some heuristics with good
performance~\cite{luNanophotonicComputationalDesign2013}.

There has been recent interest in understanding global properties of solutions for physical design problems: lower bounds for optimal design objectives in photonic design have been studied via the use of convex relaxations found by physical arguments~\cite{millerFundamentalLimitsOptical2016, shimFundamentalLimitsNearField2019}, duality theory~\cite{angerisComputationalBoundsPhotonic2019, gustafssonUpperBoundsAbsorption2019, moleskyTOperatorLimitsElectromagnetic2020}, among others~\cite{moleskyBoundsAbsorptionThermal2019}. We instead analyze a \emph{convex restriction}~\cite[\S2.1]{diamondGeneralSystemHeuristic2018} of the physical design problem, potentially providing another approach for analyzing properties of global solutions and for creating fast heuristics.


In this paper, we consider a simple (but very general) formulation of a class of physical design problems which includes problems in thermal design, photonic inverse design with scalar fields and convex objectives, and some types of control problems. This formulation offers some insights into the properties of global solutions for these problems. For example, in many practical cases, problems with linear objectives can be shown to have optimal extremal designs (in the case of physical design) or bang-bang controls (in the case of control). As another example, we observe that it suffices to know only the sign of a subset of variables in order to globally solve the problem efficiently, even though the original problem is NP-hard. The formulation also suggests a heuristic which appears to have good performance for many kinds of physical design problems, and we give numerical examples of this heuristic applied to a few different problems.

%
%
%

\section{General problem formulation}
We consider a problem of the form
\begin{equation}\label{eq:main}
\begin{aligned}
& \text{minimize} & &  f(x, u, v)\\
& \text{subject to} & & (x, u, v) \in \C\\
&&& u = \diag(\theta) v\\
&&& \tmin \le \theta \le \tmax,
\end{aligned}
\end{equation}
where $f: \reals^n \times \reals^m \times \reals^m \to \reals$ is a convex function over our variables $x \in \reals^m$ and $u, v \in \reals^n$, $\C\subseteq \reals^n \times \reals^m \times \reals^m$ is a convex constraint set, and $\theta \in \reals^n$ is our design variable whose limits are $\tmin, \tmax \in \reals^n$. While apparently simple, many physical design problems can be expressed as instances of problem~\eqref{eq:main}; we show a few examples in~\S\ref{sec:applications}. We call $(x, u, v)$ the \emph{field} (corresponding to, \eg, the electric field in photonic design) and $\theta$ the \emph{design parameters} (corresponding to, \eg, the permittivity in photonic design).  
We say that $\theta$ is \emph{extremal} whenever $\theta_i \in \{\tmin_i, \tmax_i\}$ for each $i=1, \dots, m$. The physics of the problem is encoded in the constraints $(x, u, v) \in \C$ and $u = \diag(\theta)v$.

In this problem, the convex set $\C$ can be any convex set specifying
constraints on the variables $(x, u, v)$, such as linear equality constraints.
On the other hand, the design parameters $\theta$
enter in a very specific way: as a diagonal term relating $u$ and $v$. 
Another way to say this is that each design parameter $\theta_i$
is the ratio of two field parameters, $u_i$ and $v_i$.

We note
that the problem~(\ref{eq:main})
is convex in $(x, u, v)$ whenever $\theta$ is fixed, and
convex in $(x, u, \theta)$ whenever $v$ is fixed. In practice, there has been
great success in applying heuristics for approximately minimizing instances
of~\eqref{eq:main} using this
observation~\cite{luInverseDesignNanophotonic2010}.

\paragraph{Absolute-upper-bound formulation.}
%
Problem~\eqref{eq:main} is equivalent to
\begin{equation}\label{eq:aub-main}
\begin{aligned}
& \text{minimize} & &  f(x, u, v)\\
& \text{subject to} & & (x, u, v) \in \C\\
&&& u = \diag(\bar \theta)v + \diag(\rho) w\\
&&& |w| \le |v|,
\end{aligned}
\end{equation}
where the absolute value is taken elementwise.
The variables of problem~\eqref{eq:aub-main} are $x \in \reals^m$ and $u, v, w \in \reals^n$, while $\bar\theta = (\tmax + \tmin)/2$ and $\rho = (\tmax - \tmin)/2$ are constants.
Note that $\bar\theta$ is the middle value of the physical parameter interval,
and $\rho$ is the radius, \ie, half the range or width of the interval.

The equivalence between problems~\eqref{eq:main} and~\eqref{eq:aub-main} can be seen by noting that, for every feasible $(x, u, v, w)$ for problem~\eqref{eq:aub-main} we can set,
\begin{equation}\label{eq:theta-map}
\theta_i = \begin{cases}
	\bar\theta_i + \rho_i w_i/v_i & v_i \ne 0,\\
	\bar\theta_i & \text{otherwise},
\end{cases}
\end{equation}
for $i=1, \dots, m$. Then, $(x, u, v, \theta)$ is feasible for~\eqref{eq:main}, with the same objective value. Note that, if $v_i = 0$, any choice of $\theta_i \in [\tmin_i, \tmax_i]$ would suffice.

Similarly, for any $(x, u, v, \theta)$ that is feasible for~\eqref{eq:main}, we can set
\[
w_i = \left(\frac{\theta_i - \bar\theta_i}{\rho_i}\right)v_i, \quad i=1, \dots, m,
\]
and then $(x, u, v, w)$ is feasible for problem~\eqref{eq:aub-main} with the same objective value.

We will refer to problem~\eqref{eq:aub-main} as the absolute-upper-bound formulation of problem~\eqref{eq:main}. This problem, like problem~\eqref{eq:main}, is nonconvex due to the inequality $|w| \le |v|$, and is hard to solve exactly.

\paragraph{NP-hardness.} We can reduce any mixed-integer convex program (MICP) to an instance of~\eqref{eq:aub-main}, implying that this problem is hard, as any instance of an NP-complete problem is easily reducible to instances of the MICP problem~\cite{karpReducibilityCombinatorialProblems1972}.

The reduction follows since we can force $v$ to be binary in problem~\eqref{eq:aub-main}. First, choose $\bar\theta = 0$, $\rho = \ones$ (and therefore $u=w$), and add $u = \ones$ to the constraint set. This immediately implies that $\ones \le |v|$. Adding the convex constraint $|v| \le \ones$ to the constraint set $\C$, yields $v \in \{\pm 1\}^n$, as required. Since $\C$ and $f$ can be otherwise freely chosen, the result follows.
 
\paragraph{Known signs.} If the signs of an optimal $v^\star$ are known for problem~\eqref{eq:aub-main}, then the problem becomes convex. We can see this as follows. If $s = \sign(v^\star) \in \{\pm 1\}^m$ is known, then we can solve the following convex problem~\cite[\S4]{cvxbook}:
\begin{equation}\label{eq:aub-cvx}
\begin{aligned}
& \text{minimize} & &  f(x, u, v)\\
& \text{subject to} & & (x, u, v) \in \C\\
&&& u = \diag(\bar \theta)v + \diag(\rho) w\\
&&& |w| \le s\circ v,
\end{aligned}
\end{equation}
where $s\circ v$ is the elementwise product of $s$ and $v$. Note that $v^\star$ (and its associated values of $x^\star$, $u^\star$, and $w^\star$) are feasible for this instance of~\eqref{eq:aub-cvx} since $|v^\star| = s\circ v^\star$, which implies that a solution of this instance of~\eqref{eq:aub-cvx} must be globally optimal for~\eqref{eq:aub-main}.

\paragraph{Global solution.} Note that problem~\eqref{eq:aub-cvx} generates a family of optimization problems over the set of possible signs, $s \in \{\pm 1\}^m$. This suggests a simple, if inefficient, way to globally solve problem~\eqref{eq:aub-main} and therefore problem~\eqref{eq:main}: solve problem~\eqref{eq:aub-cvx} for the $2^m$ possible signs, $s \in \{\pm 1\}^m$, to obtain optimal values $p^\star(s)$ for each set of signs $s$. A solution $(x^\star, u^\star, v^\star, w^\star)$ for any optimal set of signs, $s^\star \in \argmin_{s \in \{\pm 1\}^m} p^\star(s)$, is then a solution to~\eqref{eq:aub-main} and therefore to~\eqref{eq:main}.

Of course, this algorithm may not be useful in practice for anything but the smallest values of $m$, but it implies that solving problem~\eqref{eq:main} requires solving only a finite number of convex problems.

\paragraph{Extremality principle.}\label{sec:maximality-principle} The rewriting given in~\eqref{eq:aub-cvx} also yields an interesting insight. If problem~\eqref{eq:aub-cvx} is a feasible linear program and $\C$ is an affine set with $\{u \mid (x, u, v) \in \C\} = \reals^m$, \ie, for each $u \in \reals^m$ there exists a $v \in \reals^m$ and an $x \in \reals^n$ such that $(x, u, v) \in \C$, then there exists a solution of~\eqref{eq:aub-cvx} such that all entries of the inequality $|w| \le s\circ v$ hold at equality. (See, \eg,~\cite[\S2.6]{bertsimasIntroductionLinearOptimization1997}.) This rewriting then implies that there exists an optimal design for which $\theta$ is extremal, by~\eqref{eq:theta-map}. A numerical example of this principle is found in~\S\ref{sec:numerical-static-diffusion}.

\section{Sign flip descent}
\label{sec:sign-flip-descent}
Since problem~\eqref{eq:aub-cvx} generates a family of optimization problems parametrized by the sign vector $s \in \{\pm 1\}^m$, we can view the original physical design problem~\eqref{eq:main} as a problem of choosing an optimal Boolean vector. A simple way of approximately optimizing~\eqref{eq:aub-main} is: at each iteration $i$, start with some sign vector $s^i \in \{\pm 1\}^m$ and solve~\eqref{eq:aub-cvx} to obtain an optimal value $p^i$. We then consider a rule for proposing a new sign vector, say $\tilde s^i \in \{\pm1\}^m$, for which we again solve~\eqref{eq:aub-cvx} and then obtain a new optimal value $\tilde p^i$. If $\tilde p^i < p^i$, we then keep this new sign vector, \ie, we set $s^{i+1} = \tilde s^i$, and repeat the procedure; otherwise, we discard $\tilde s^i$ by setting $s^{i+1} = s^i$, and repeat the procedure, proposing a new sign vector in the next iteration. This is outlined in algorithm~\ref{alg:main}.

\begin{algorithm}[H]
\caption{Sign flip descent}
\label{alg:main}
\begin{algorithmic}
\State Start with a feasible initial sign vector $s^1 \in \{\pm 1\}^m$.
\State {\emph{Optimize.} Solve problem~\eqref{eq:aub-cvx} with signs $s^1$ to receive objective value $p^1$.}
\For {$k = 1, \dots, n_\mathrm{iter}$}
\State {\emph{Propose.} Propose a new set of signs $\tilde s^k \in \{\pm 1\}^m$.}
\State {\emph{Optimize.} Solve~\eqref{eq:aub-cvx} with the sign vector $\tilde s^k$ to receive objective value $\tilde p^k$.}
\If {$\tilde p^k < p^k$}
  \State {$s^{k+1} = \tilde s^k$.}
\Else
  \State {$s^{k+1} = s^k$.}
\EndIf
\EndFor
\State \textbf{return} $s^{n_\textrm{iter}}$.
\end{algorithmic}
\end{algorithm}

By construction, any algorithm of the form of algorithm~\ref{alg:main} is a descent algorithm since each iteration is feasible and the objective value is decreasing on each iteration. We outline two possible rules for proposing new sets of signs at each iteration.

\paragraph{Greedy sign rule.}
A simple rule for choosing signs is to begin at iteration $k$ with some set of signs $s^k$. We then define a new set of signs $\tilde s^k$ with $\tilde s^k = s^k$ except at the $k$th entry where we have $\tilde s^k_k = -s^k_k$ (or, if $k > m$ then we pick the entry at index $1 + (k-1\mod m)$, \ie, such that the entries are changed, one-by-one, in a round-robin fashion). We stop whenever flipping any one entry of $s^k$ does not yield a lower objective value.

The greedy sign rule has two useful properties. First, the rule guarantees local optimality in the following sense: if algorithm~\ref{alg:main} returns $s^\star$, then changing any one sign of $s^\star$ will not decrease the objective value. Second, the rule terminates in finite time, since the corresponding algorithm is a descent algorithm and there are a finite number of possible sign vectors. On the other hand, the algorithm is often slow for anything but the smallest designs: to reach the terminating condition, we have to solve at least $m$ convex optimization problems.

\paragraph{Field-based rule.}
Another simple rule that appears to work very well in practice is based on the observation that, for many choices of sign vectors $s^k$, the inequality $|w| \le s^k \circ v$ has many entries of $v$ that are zero. If $v_i$ is zero for some index $i=1, \dots, m$, this suggests that the sign $s_i^k$ might have been originally set incorrectly: in this case, we propose a new vector $\tilde s^k$ which is equal $s^k$ at all entries $i=1, \dots, m$ for which $v_i$ is nonzero and has opposite sign at all entries $i$ for which $v_i$ is zero.

Note that this new proposed vector will always have an optimal value $\tilde p^k$ which is at least as small as the optimal value for $s^k$, \ie, $\tilde p^k \le p^k$. This observation, coupled with the proposed rule, suggests that we should stop whenever there are no signs left to flip, or whenever the iterations stop decreasing as quickly as desired, \ie, whenever $p^k - p^{k+1} < \eps$.

While this rule does not necessarily guarantee local optimality, it always terminates in finite time with the given stopping conditions and appears to work well in practice (requiring, in comparison to the greedy sign rule, much fewer than $m$ iterations to terminate) as shown in~\S\ref{sec:numerical-examples}.

\section{Applications}
\label{sec:applications}
We describe a few interesting design problems that reduce to problems of the form of~\eqref{eq:main}.

\subsection{Diagonal physical design}
\label{sec:diagonal-physical-design}
As in, \eg,~\cite{angerisComputationalBoundsPhotonic2019}, many physical design problems can be written in the following way: 
\begin{equation}\label{eq:physical-design}
\begin{aligned}
& \text{minimize} & &  f(z)\\
& \text{subject to} & & (A + \diag(\theta))z = b\\
&&& \theta^\mathrm{min} \le \theta \le \theta^\mathrm{max},
\end{aligned}
\end{equation}
where $A \in \reals^{n\times n}$ describes the physics of the problem, while $b \in \reals^n$ describes the excitation, and $\theta \in \reals^n$ are the design parameters of the system, chosen to minimize some convex objective function $f: \reals^n \to \reals$ of the field $z \in \reals^n$. Our variables in this problem are the field $z$ and the design parameters $\theta$.

We can write a problem of the form of~\eqref{eq:physical-design} as a problem of the form~\eqref{eq:main} by introducing a new variable $u$ with constraint $u = \diag(\theta)z$ and rewriting the equality constraint of~\eqref{eq:physical-design} with this new variable, $Az + u = b$. As the set of $(z, u)$ satisfying $Az + u = b$ forms a convex (in fact, affine) set, the resulting problem,
\[
\begin{aligned}
& \text{minimize} & &  f(z)\\
& \text{subject to} & & Az + u = b\\
&&& u = \diag(\theta)z\\
&&& \theta^\mathrm{min} \le \theta \le \theta^\mathrm{max},
\end{aligned}
\]
is of the form of~\eqref{eq:main} which can be easily rewritten into the form of~\eqref{eq:aub-main}.

\subsection{Static diffusion design}
\label{sec:static-diffusion}
Consider a flow problem on a graph $G = (V, E)$ where we choose the conductance $g_k \in \reals$ across each edge $k \in E$, constrained to satisfy $\gmin_k \le g_k \le \gmax_k$, to minimize some function $f: \reals^{|V|} \to \reals$ of the potentials $e \in \reals^{|V|}$, given some sources $s \in \reals^{|V|}$.

To compactly write the conditions this system must satisfy, let the matrix $A \in \reals^{|V|\times |E|}$ be the incidence matrix for the graph $G$ defined to be (see~\cite[\S7.3]{linalg}):
\[
A_{ij} = \begin{cases}
  +1 & \text{edge $j$ points to node $i$}\\
  -1 & \text{edge $j$ points from node $i$}\\
  0 & \text{otherwise}.
\end{cases}
\]
We can then write the steady-state diffusion equation as
\begin{equation}\label{eq:heat-steady}
A\diag(g)A^Te = s,
\end{equation}
where $A\diag(g)A^T$ can be recognized as the graph Laplacian of $G$ with edge weights $g$. This equation can also be seen as the discrete form of the heat equation 
on a graph $G$~\cite{solomonPDEApproachesGraph2015}.

The corresponding optimization problem is then an instance of~\eqref{eq:main}:
\begin{equation}\label{eq:static-diffusion}
\begin{aligned}
& \text{minimize} & &  f(e)\\
& \text{subject to} & & v = A^Te\\
&&& Aw = s\\
&&& w = \diag(g)v\\
&&& \gmin \le g \le \gmax,
\end{aligned}
\end{equation}
where we have introduced two new variables $w, v \in \reals^{|E|}$, in addition to the potential $e \in \reals^{|V|}$ and the conductances $g \in \reals^{|E|}$. As before, $A \in \reals^{|V|\times |E|}$ is the incidence matrix, $s \in \reals^{|V|}$ are the sources at each node, while $c \in \reals^{n}$ is a vector such that $c^Te$ is the average temperature over the desired region.

\subsection{Dynamic diffusion control}
\label{sec:temperature-control}
Similarly to~\S\ref{sec:static-diffusion}, we can consider the time-varying generalization of~\eqref{eq:heat-steady} given by
\[
Ce_{t+1} = Ce_t - h A\diag(g_t) A^Te_t + hBu_t,
\]
at each time $t=1, \dots, T$, with step size $h > 0$. Here, $c \in \reals^{|V|}_{++}$ is the heat capacity of each node and $C = \diag(c)$, while $u_t \in \reals^{n}$ are the inputs given to the system, $B \in \reals^{|V| \times n}$ is a matrix mapping these inputs to the power added or removed from each node, $g_t \in \reals^{|V|}$ are the conductances at each node, and $e_t \in \reals^{|V|}$ is the temperature at each node.

In this case, we can minimize any convex function of the temperatures and inputs by appropriately choosing the conductances and inputs:
\begin{equation}\label{eq:dynamic-control}
\begin{aligned}
& \text{minimize} & &  f(e, u)\\
& \text{subject to} & & Ce_{t+1} = Ce_t + hAw_t + hBu_t, \quad t=1, \dots, T-1\\
&&& v_t = A^Te_t,\quad t=1, \dots, T-1\\
&&& w_t = \diag(g_t)v_t, \quad t=1, \dots, T-1\\
&&& \gmin \le g_t \le \gmax, \quad t=1, \dots, T-1,
\end{aligned}
\end{equation}
where as before, we have introduced the variables $v_t, w_t \in \reals^{|E|}$, for each $t=1, \dots, T$.

We can see problem~\eqref{eq:dynamic-control} as a nontraditional control problem. A particular example is: we have a set of rooms with temperatures $e_t$ at time $t$ which we wish to keep within some comfortable temperature range. We are allowed to open and close vents (equivalently, change the conductances $g_t$ at each time $t$) and turn on and off heat pumps (via the control variable $u_t$), while paying a cost for the latter. A simple question could be: what is an optimal set of actions such that the input cost is minimized while keeping the temperatures $e_t$ within some specified bounds? We show a simple example of this in~\S\ref{sec:numerical-dynamic-control}.

\section{Numerical examples}
\label{sec:numerical-examples}
Julia~\cite{bezansonJuliaFreshApproach2017} code for all examples in this section is available in the following Github repository: \texttt{https://github.com/angeris/pd-heuristic}. We use the JuMP modeling language~\cite{dunningJuMPModelingLanguage2017} to interface with Mosek~\cite{mosek}. All times reported are on a 2015 2.9~GHz dual-core MacBook Pro.

\subsection{Photonic design}
In this example, we wish to choose the speed of a wave satisfying Helmholtz's equation at each point in some domain $\Omega \subseteq \reals^2$ in order to minimize a convex function of the field.

\paragraph{Helmholtz's equation.} More specifically, the speed of the wave $c : \Omega \to \reals_{++}$ is chosen such that the field $\psi : \Omega \to \reals$ at a specific frequency $\omega \in \reals_+$ with excitation $\phi : \Omega \to \reals$ satisfies Helmholtz's equation,
\begin{equation}\label{eq:helmholtz}
	\nabla^2 \psi(x, y) + \left(\frac{\omega}{c(x, y)}\right)^2\psi(x, y) = \phi(x, y),
\end{equation}
at each point $(x, y) \in \Omega$. Additionally, we require that the chosen speeds are bounded such that $0 < c^\mathrm{min}(x, y) \le c(x, y) \le c^\mathrm{max}(x, y)$ at each point $(x, y) \in \Omega$, and we assume Dirichlet boundary conditions such that $\psi(x, y) = 0$ for $(x, y) \in \partial \Omega$, \ie, we require the field to be zero at every point on the boundary of the domain.

In this case (as in~\cite[\S5.1]{angerisComputationalBoundsPhotonic2019}), we will work with a discretized form of~\eqref{eq:helmholtz} where $z \in \reals^n$ is the discretized field ($\psi$), $b \in \reals^n$ is the discretized excitation ($\phi$), $\theta \in \reals^n$ is the discretized speed of the wave ($c$), and $A \in \reals^{n\times n}$ is the discretized version of the Laplacian operator ($\nabla^2$), such that
\begin{equation}\label{eq:discretized}
Az + \diag(\theta)z = b,
\end{equation}
approximates~\eqref{eq:helmholtz} at each point $(x_i, y_i) \in \Omega$ for $i=1,\dots, n$.

\begin{figure}
    \centering
    \includegraphics[width=\textwidth]{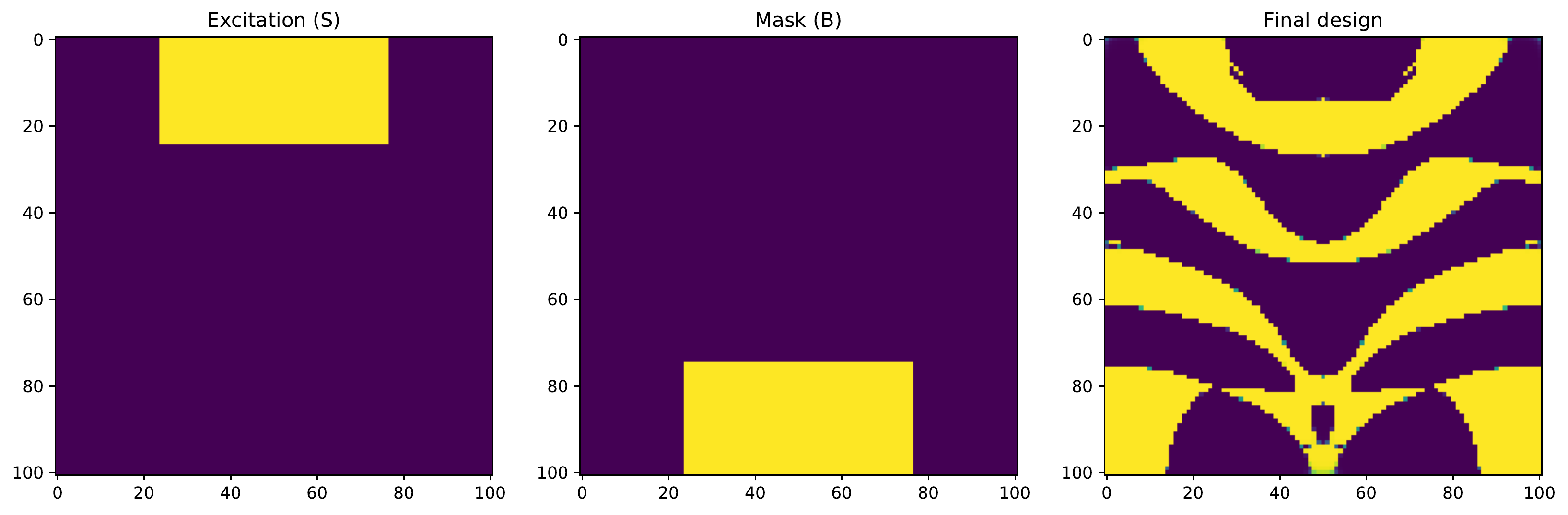}
    \caption{Approximately optimal photonic design. The leftmost figure specifies $S \subseteq \{1, \dots, n\}$ in yellow, the center specifies $B \subseteq \{1, \dots, n\}$, while the rightmost figure gives the design, $\theta$.}
    \label{fig:photonic-design}
\end{figure}

\paragraph{Problem data.} In this case, the problem data are given by $w = 4\pi$, with $n = 101 \times 101 = 10201$, while the convex objective function $f: \reals^n \to \reals$ is given by
\[
f(z) = \sum_{i \in B} z_i^2,
\]
where $B \subseteq \{1, \dots, n\}$ is the box indicated in figure~\ref{fig:photonic-design}, and the excitation $b$ is defined as
\[
b_i = \begin{cases}
	1 & i \in S\\
	0 & \text{otherwise},
\end{cases}
\]
for each $i=1, \dots, n$, where $S \subseteq \{1, \dots, n\}$ is the box indicated in figure~\ref{fig:photonic-design}. Here, $\tmin = 1$ and $\tmax = 2$. We set the tolerance parameter of the algorithm to $\eps = 10^{-4}$. We initialize the algorithm by finding a solution to equation~\eqref{eq:discretized} with $\theta = (\tmax + \tmin)/2$ and use the signs of this solution as the initial sign vector.

\paragraph{Numerical results.} With the given problem data, the algorithm terminates at 102 iterations with a total time of about 4 minutes, roughly around 2 seconds per iteration. This time could be very much shortened since the current implementation does not warm-start any of the current iterations, essentially solving the problem from scratch at each iteration. The final design is shown in figure~\ref{fig:photonic-design}.

\subsection{Thermal design}
\label{sec:numerical-static-diffusion}
In this design problem, as in~\S\ref{sec:static-diffusion}, we seek to set the conductances on a graph in order to minimize the average temperature of a subset of points in the center of a 2D grid of size $m\times m$, given a heat source and a heat sink at opposite corners of the 2D grid. This is an instance of the diffusion problem where $A \in \reals^{|E| \times |V|}$ is the incidence matrix of the grid and $s \in \reals^{|V|}$ are the heat sources and sinks. This problem can be written as an instance of~\eqref{eq:static-diffusion} where the potentials $e \in \reals^{|V|}$ are the temperatures at each point in the grid.

\begin{figure}
    \centering
    \includegraphics[width=\textwidth]{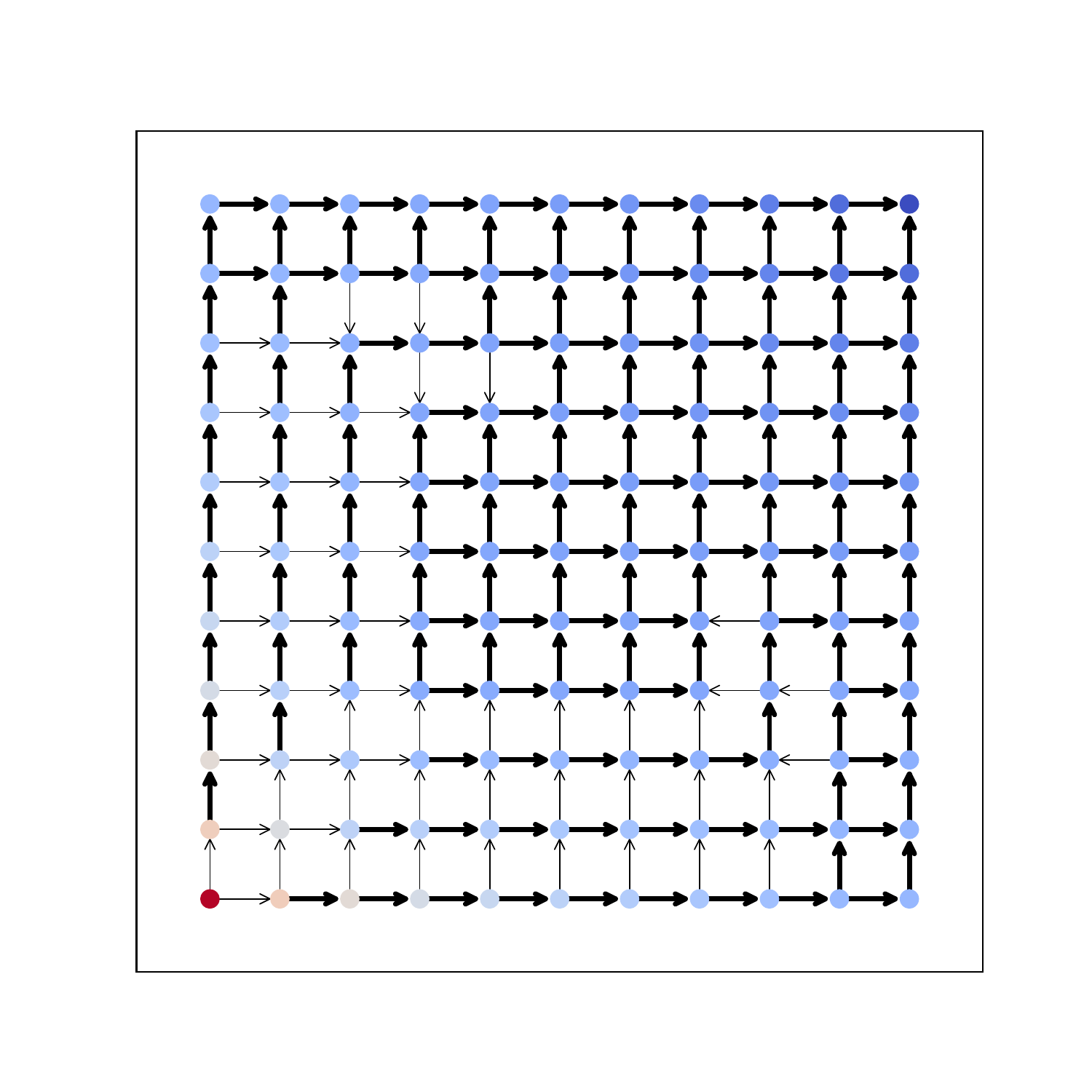}
    \caption{Approximately optimal design for $m=11$. Arrows indicate the direction of flow used for this design, colors indicate the temperature at each node, while edge thickness indicates the conductance at each edge.}
    \label{fig:small-design}
\end{figure}

\begin{figure}
    \centering
    \includegraphics[width=\textwidth]{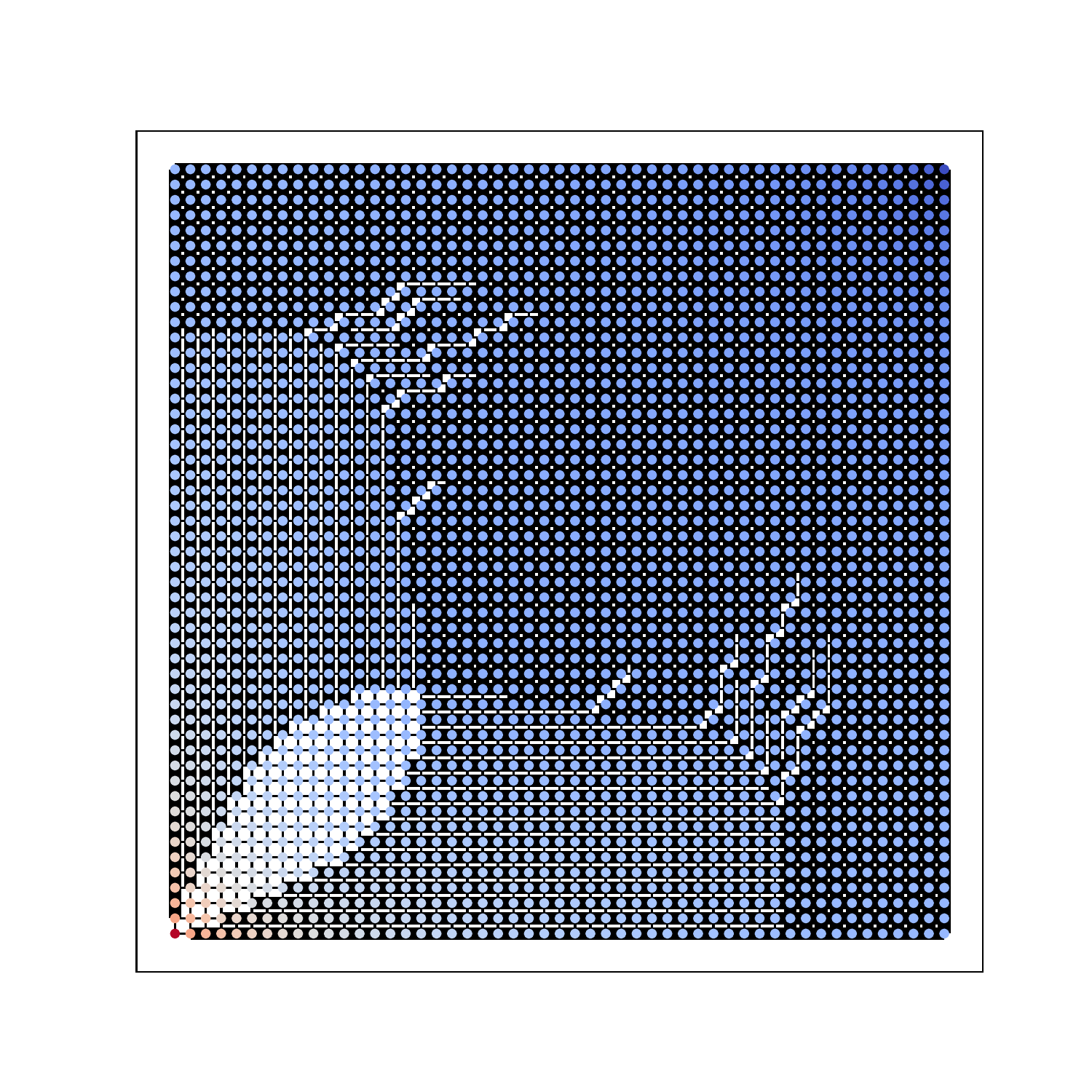}
    \caption{Approximately optimal design for $m=51$.}
    \label{fig:large-design}
\end{figure}

\paragraph{Problem data.}
Our convex objective function $f: \reals^{|V|} \to \reals$ is given by
\[
f(e) = c^Te,
\]
where $c \in \reals^{|V|}$ is a vector such that $c_i = 1$ if vertex $i$ lies in the center square of size $\lfloor(m-1)/4\rfloor \times \lfloor(m-1)/4\rfloor$ (indicated as a grey box in figures~\ref{fig:small-design} and~\ref{fig:large-design}) while $c_i = 0$ otherwise. There is a heat source set at the bottom left corner of the grid and a heat sink set at the top right corner of the grid. We set the minimal and maximal conductances as $\gmin = 1$ and $\gmax = 10$ at each edge.

We approximately optimize the conductances in this problem by using the field-based heuristic described in~\S\ref{sec:sign-flip-descent}. The directions are initialized by solving the problem with uniform conductances.

\paragraph{Numerical results.} A small example is given in figure~\ref{fig:small-design} with $m=11$ (which shows the chosen directions of flow), while a relatively large design is given in figure~\ref{fig:large-design} with $m=51$. In both figures, thick edges indicate that conductance is maximized at that edge while thin edges indicate that conductance is minimized (see the extremality principle in~\S\ref{sec:maximality-principle} for more details). The color of each node indicates the potential value, with red values indicating a higher potential and blue values indicating a lower one. We note that our heuristic recovers similar tendril-like patterns to those found in, \eg,~\cite[\S4]{haertelDesignThermalSystems2018}.

With the provided data, the heuristic terminates after 7 iterations, taking a total time of around .4 seconds in the case with $m=11$, with an objective value of about .115. The case with $m=51$ terminates after 14 iterations, taking a total time of around 20.5 seconds with an objective value of approximately .239.

\subsection{Temperature control}\label{sec:numerical-dynamic-control}
In this example, we wish to keep the temperature of two rooms in a range of desired temperatures by appropriately closing and opening vents to the outside and between rooms and turning heat pumps on and off at specified times, while minimizing the total power consumption. We will also require that the controls and the temperatures be periodic.

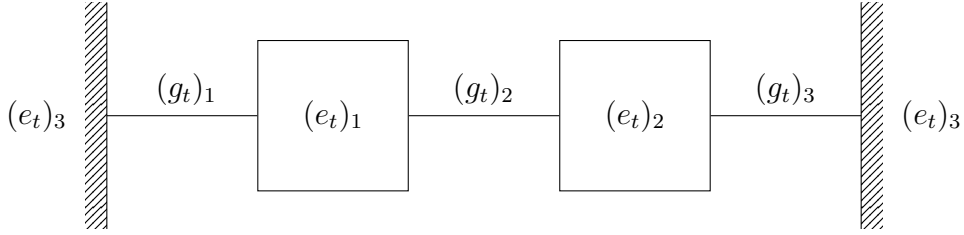
\begin{figure}
    \centering
    \begin{tikzpicture}
    	\tikzstyle{ground}=[fill,pattern=north east lines,draw=none,minimum width=0.5,minimum height=3cm]
    	\tikzstyle{room}=[draw,minimum height=2cm,minimum width=2cm]
    	
    	\node (left_label){$(e_t)_3$};
    	\node (temp_ambient)[ground,right=.1 of left_label]{};
    	\draw (temp_ambient.north east) -- (temp_ambient.south east);
    	
    	\node (room1)[room,right = 2 of temp_ambient] {$(e_t)_1$};
    	\draw (temp_ambient) -- ++(1.2, 0) node[above] {$(g_t)_1$} -- (room1);
    	
    	\node (room2)[room,right = 2 of room1] {$(e_t)_2$};
    	\draw (room1) -- ++(2, 0) node[above] {$(g_t)_2$} -- (room2);
    	
    	\node (temp_ambient_right)[ground, right = 2 of room2] {};
    	\draw (temp_ambient_right.north west) -- (temp_ambient_right.south west);
    	\draw (room2) -- ++(2, 0) node[above] {$(g_t)_3$} -- (temp_ambient_right);
    	
    	\node [right=.1 of temp_ambient_right] {$(e_t)_3$};
    \end{tikzpicture}
    \caption{Graph set up for the temperature control problem. Here, $(e_t)_3$ is the ambient temperature at time $t$, while $(e_t)_1$ and $(e_t)_2$ are the temperatures of rooms 1 and 2, respectively. The $g_t$ are the conductances of the indicated edges.}
    \label{fig:room-graph}
\end{figure}

\paragraph{Problem data.} We can write this as an instance of problem~\eqref{eq:dynamic-control} with
\[
B = .2 I, \quad C = \diag((.3, .1)), \quad \gmin = 1, \quad \gmax = 10,
\]
and $A$ is the incidence matrix of the graph shown in figure~\ref{fig:room-graph}, while
\[
(e_t)_3 = 70 + 20\sin\left(\frac{4\pi t}{T}\right), \quad t=1, \dots, T,
\]
where $T = 300$. Since we will require that the room temperatures be periodic, we then have
\[
(e_1)_1 = (e_T)_1, \quad (e_1)_2 = (e_T)_2.
\]
Finally, we will require that the temperatures remain in some a range,
\[
65 \le (e_t)_1, (e_t)_2 \le 75, \quad \quad t = 1, \dots, T,
\]
while minimizing
\begin{equation}\label{eq:obj-control}
f(e, u) = h\|u\|_2 + \eta h \sum_{t=1}^{T-1}\|e_{t+1} - e_t\|_2,
\end{equation}
where $h=1/T$ and $\eta = 10^{-4}$ is a small regularization parameter that ensures the resulting trajectories are smooth.

We initialize the problem with the signs given by assuming that $g_t = (\gmin + \gmax)/2$ for all $t=1, \dots, T-1$ and using the heat pumps $u_t$ to ensure the temperature in the rooms remains above 65 and below 75.

\begin{figure}
    \centering
    \includegraphics[height=.9\textheight]{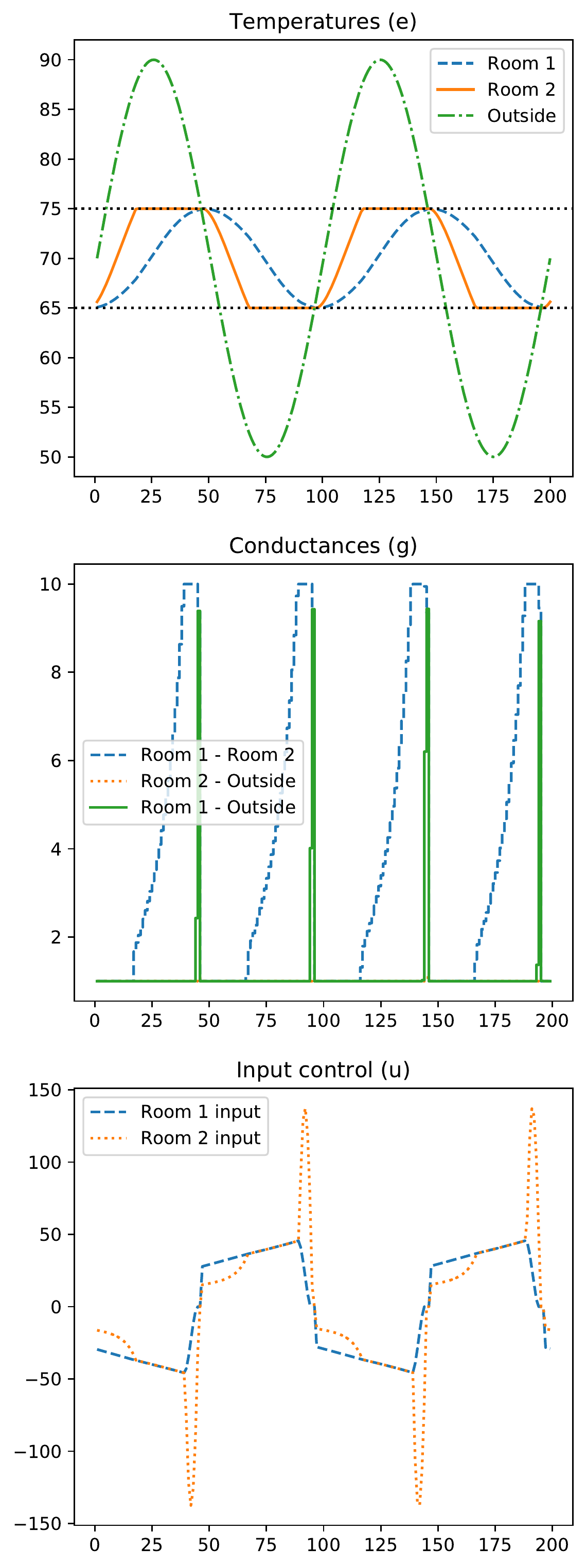}
    \caption{Approximately optimal control.}
    \label{fig:control}
\end{figure}

\paragraph{Numerical results.} We approximately optimize this instance using the field-based heuristic of~\S\ref{sec:sign-flip-descent}, with the result shown in figure~\ref{fig:control}. With the provided data, the heuristic terminates in 3 iterations, with a total time of around 1.56 seconds. The final approximately optimized problem has an objective value of around $836$.

\section{Conclusion}
This paper presented a new problem formulation and an associated heuristic which may be of practical use for a general class of physical design problems, which appears to have good practical performance on many different kinds of physical design problems. Additionally, this problem formulation implies a few interesting facts, most notably that the class of problems can be efficiently solved even when only the signs of an optimal solution are known and that, in a few important cases, there exist globally optimal maximal designs.

\paragraph{Future work.} There are several notable exceptions to the class of problems which are included in the formulation given in~\eqref{eq:main}, with the most important being designs whose parameters are constrained to be equal. This means that, at the moment, a direct application to photonic design in three dimensions, circuit design with complex impedances, or multi-scenario physical design, is not possible with the current problem formulation. We suspect a suitable generalization of~\eqref{eq:main} might yield similarly interesting insights and, potentially, new heuristics for physical design.

\section*{Acknowledgements}
The authors would like to acknowledge Akshay Agrawal, Shane Barratt, and Rahul Trivedi for helpful comments and edits.

\appendix

\bibliography{citations.bib}

\end{document}